\documentclass[12pt]{article}
\usepackage{amsmath}
\usepackage{amssymb}
\usepackage[cp1251]{inputenc}
\usepackage[russian]{babel}
\newcommand{\liml}{\lim\limits_{t \to + \infty}}
\newcommand{\supl}{\sup\limits_{t\ge t_0}}
\newcommand{\Supl}[1]{\sup\limits_{#1}}
\newcommand{\expb}[1]{\exp\biggl\{#1\biggr\}}
\newcommand{\il}[2]{\int\limits_{#1}^{#2}}
\newcommand{\ilp}[1]{\int\limits_{#1}^{+\infty}}

\newcommand{\ph}{\phantom{a}}
\newcommand{\phh}{\phantom{aaa}}

\newcommand{\re}{Re \hskip 1.5pt}
\newcommand{\sist}[2]{\left\{
\begin{array}{l}
{#1}\\
\ph\\
{#2}
\end{array}
\right.}

\begin{document}

MSC 34D20

\vskip 20pt

\centerline{\bf  Stability criterion for linear systems}
\centerline{\bf  of  ordinary differential equations}

\vskip 10 pt

\centerline{\bf G. A. Grigorian}

\vskip 10 pt

\centerline{0019 Armenia c. Yerevan, str. M. Bagramian 24/5}
\centerline{Institute of Mathematics NAS of Armenia}
\centerline{E - mail: mathphys2@instmath.sci.am, \ph phone: 098 62 03 05, \ph 010 35 48 61}

\vskip 20 pt

\noindent
Abstract. The Riccati equation method is used to establish a new stability  criteria for linear systems of  ordinary differential equations. Two examples are presented in which the obtained result is compared with the results obtained by the Lyapunov and Bogdanov methods, by a method involving estimates of solutions in the Lozinskii logarithmic norms and by the freezing method.

\vskip 20 pt


\vskip 20 pt
{\bf  1. Introduction}. Denote by $\mathbb{H}$ the algebra of quaternions $q = q_0 + i q_1 + j q_2 + k q_3$, where $q_n\in \mathbb{R}, \ph n=\overline{1,3}, \ph i, \ph, j, \ph k$ are the imaginary units satisfying the relations
$$
i^2 = j^2  = k^2 = ijk = -1, \ph ij = -ji = k, \ph jk = -kj = i , \ph ki = - ik = j.
$$
$\mathbb{H}$ is an Euclidean space endowed with the norm (modulus) $|q| =\sqrt{q_0^2 + q_1^2 + q_2^2 + q_3^2}$ of its elements $q= q_0 + i q_1 + j q_2 + k q_3$. Then $\mathbb{H}^n \equiv \{\mathfrak{q}=(\mathfrak{q}_1, ..., \mathfrak{q}_n): \mathfrak{q}_l\in \mathbb{H}, \ph l= \overline{1,n}\}$ is a $4n$ dimensional Euclidean space endowed with the norm $||\mathfrak{q}|| =\sqrt{|\mathfrak{q}_1|^2 + ... + |\mathfrak{q}_n|^2}$ of its elements $\mathfrak{q}=(\mathfrak{q}_1, ..., \mathfrak{q}_n)$.
Let $J$ be a quaternion with $J^2 = -1$. Denote by $\mathbb{C}_J$ the set of quaternions $q$ of the form $q=\alpha + J\beta, \ph \alpha, \beta \in \mathbb{R}$ . Obviously there is an isometria  $\mathbb{C}_j\stackrel{f_J}{\rightarrow} \mathbb{C}$ between $\mathbb{C}_J$ and the field of the complex numbers $\mathbb{C}$. Denote by $\mathbb{C}_J^{n\times n}$ the set of the matrices of dimension $n\times n$ with the elements from $\mathbb{C}_J$. The isometria $f_J$ generates an isometria $F_j : \mathbb{C}_J^{n \times n} \rightarrow \mathbb{C}^{n\times n}$ between $ \mathbb{C}_J^{n \times n}$ and $ \mathbb{C}^{n \times n}$ by the following rule: if $M = (q_{lm})_{l,m =1}^n, \ph q_{lm}\in \mathbb{C}_J, \ph l,m = \overline{1,n}$, then $F_J(M) = (f_J(q_{lm}))_{l,m=1}^n$. A matrix $U\in \mathbb{C}_J^{n\times n}$ we will call $J$-unitary if $U U^* = U^* U = I$, where $U^*$ is the transpose conjugate for $U$, $I$ is the identity matrix of dimension $n\times n.$
 Let $\mathcal{A}: \mathbb{H}^m \rightarrow \mathbb{H}^n$ be a linear continuous operator acting from $\mathbb{H}^m$ to  $\mathbb{H}^n$.
 Denote by $||\mathcal{A}||$ the norm of $\mathcal{A} ; \ph ||\mathcal{A}||=\Supl{q\ne0}\frac{||\mathcal{A}\mathfrak{q}||}{||\mathfrak{q}||}$. Due to the isometria $F_J$ for any $J$-unitary matrix $U$ the equalities
$$
||U|| = || U^*|| = 1
$$
are valid.

Let $A(t), \ph B(t), \ph C(t)$ and $D(t)$ be quaternionic-valued  continuous matrix functions  of dimensions $m\times m, \ph m\times n, \ph n\times m,$ and $n\times n$ respectively on $[t_0,+\infty)$. Consider the linear system
$$
\sist{\phi' = A(t) \Phi + B(t) \Psi,} {\Psi' = C(t) \Phi + D(t) \Psi, \ph t \ge t_0.} \eqno (1.1)
$$
By a solution of this system we mean an ordered pair $(\Phi(t), \Psi(t))$ of  continuously differen-\linebreak tiable vector functions $\Phi(t)=(\phi_1(t), ..., \phi_m(t))$ and $\Psi(t)= (\psi_1(t), ..., \psi_n(t))$, satisfying (1.1)  on $[t_0,+\infty)$.

Due to the method, used in this paper we must put the following restrictions on the coefficients of the system (1.1).

\noindent
$a) \ph A(t)\il{t_0}{\tau}A(s) d s = \il{t_0}{\tau}A(s) d s A(t), \ph   \ph D(t)\il{t_0}{\tau}D(s) d s = \il{t_0}{\tau}D(s) d s D(t), \ph  \tau,  t \ge t_0;$

\noindent
the matrices $\il{\tau}{t}A(s) d s$ and $\il{\tau}{t}D(s) d s$ are respectively $J_A$-unitary and $J_D$-unitary \linebreak ($J_A^2 = J_D^2 = -1$)  equivalent to some diagonal matrices,~ i. e.

\noindent
$b) \ph U_A(\tau;t) \il{\tau}{t}A(s) d s U_A(\tau;t)^* = diag\{a_1(\tau;t), ..., a_m(\tau;t)\},$

$\phantom{aaaaaaaaaaaaaaaaaaaaaaaaaa} U_D(\tau;t) \il{\tau}{t}A(s) d s U_D^*(\tau;t) = diag\{d_1(\tau;t), ..., d_m(\tau;t)\}$,

\noindent
where $U_A(\tau;t)$ and $U_D(\tau;t)$ are some $J_A$-unitary and $J_D$-unitary  matrix functions respec-\linebreak tively, $U_A^*(\tau;t)$ and $U_D^*(\tau;t)$ are their transpose conjugate respectively, $a_l(\tau;t), \linebreak d_r(\tau;t) \in~ \mathbb{H},\ph l=\overline{1,m}, \ph r= \overline{1,n}, \ph t \ge~ \tau \ge  t_0,$

{\bf Remark 1.1.} {\it A quaternionic-valued square matrix $V$ we will call normal if it is permutable with its transpose conjugate $V^*: \ph V V^*  = V^* V$. For example the matrices $V_\pm\equiv \begin{pmatrix}a&b\\ \pm a&b\end{pmatrix}$, where $a, \ph b \in \mathbb{C}_J$ for some $J\in \mathbb{H}$ with $J^2 =-1$ are normal. The cyclic matrix  $V \equiv  \begin{pmatrix}0&1&0&\dots&0&0\\
                                               0&0&1&\dots&0&0\\
                                               .&.&.&\dots&.&.\\
                                               0&0&0&\dots&0&1\\
                                                1&0&0&\dots&0&0\end{pmatrix}$ is also  normal. Obviously any quaternionic-valued Hermitian matrix $V_H = V_H^*$ and any quaternionic-valued skew-symmetric matrix $V_S=-V_S^*$ are normal as well.
For any $J\in \mathbb{H}$ with $J^2 =-1$ and a normal matrix $V$ set $\Omega_{J,V}^{p\times p} \equiv \Bigl\{ \sum\limits_{l=0}^Nc_l(t) V^l : \ph c_l(t) \in \mathbb{C}_J, \ph l=\overline{1,p}, \ph t \ge t_0\Bigl \}.$ One can show that the conditions $a)$ and $b)$ are satisfied if in particular $A(t)\in \Omega_{J_1,V_1}^{m\times m}, \ph D(t)\in \Omega_{J_2,V_2}^{n\times n}$ for some $J_1, J_2 \in \mathbb{H}$ with $J_1^2 = J_2^2 = -1$ and normal matrices $V_1 \in \mathbb{C}_{J_1}^{m\times m}$ and $V_2 \in \mathbb{C}_{J_2}^{n\times n}$.}

{\bf Definition 1.1}. {\it The system (1.1) is called Lyapunov (asymptotically) stable if its all solutions are bounded on $[t_0,+\infty)$ (vanish at $+\infty$).}

Study of the question of stability  of the system (1.1), in general, of linear systems of ordinary differential equations, is an important problem of qualitative theory of differential equations. Being of interest not only in theory but also for applications it is the subject of numerous investigations (see e. g., [1 - 12]). There exist many methods of estimation of solutions of linear systems of ordinary differential equations allowing to describe (to detect) wide classes of stable and (or) unstable systems of ordinary differentia; equations. Among them the main ones include the Lyapunov's, Bogdanov's, Lozinski's estimate methods and the freezing method (see [4], pp. 40 -98, 132 -145). The fundamental method of Lyapunov characteristic exponents allows  to describe the asymptotic growth of solutions of linear systems of ordinary differential equations via these exponents and therefore carrying out the stability behavior of solutions of the system. However the application of this method has some difficulties, arising in the calculation process of Lyapunov characte-\linebreak ristic exponents. There exist also other estimation methods for special classes of linear systems of ordinary differential equations (see e. g., [5 - 10]), allowing to describe wide classes of stable and (or) unstable linear systems of ordinary differential equations. Hovever these, indicated above and other methods cannot completely describe the  stable and unstable linear systems of ordinary differential equations (in terms of their coefficients).

In this paper we use the Riccati equation method to establish a new stability criterion for the system (1.1). By two examples we compare the obtained result with the results obtained by the Lyapunov and Bogdanov methods, by a method involving estimates of solutions in the Lozinskii logarithmic norms and by the freezing method.

{\bf 2. Auxiliary propositions}.
 Let $f(t), \ph g(t), \ph h(t), \ph f_1(t), \ph g_1(t), \ph h_1(t)$ be real-valued  continuous  functions on $[t_0,+\infty)$. Consider the Riccati equations
$$
y' + f(t) y^2 + g(t) y + h(t) = 0, \phh t \ge t_0. \eqno (2.1)
$$
$$
y' + f_1(t) y^2 + g_1(t) y + h_1(t) = 0, \phh t \ge t_0, \eqno (2.2)
$$
$j=1,2,$ and the differential inequalities
$$
\eta' + f(t) \eta^2 + g(t) \eta + h(t) \ge 0, \phh t \ge t_0. \eqno (2.3)
$$
$$
\eta' + f_1(t) \eta^2 + g_1(t) \eta + h_1(t) \ge 0, \phh t \ge t_0, \eqno (2.4)
$$
Note that for $f(t) \ge 0 \phantom{a} (f_1(t) \ge 0), \phantom{a} t\ge t_0$  every solution of the linear equation  $\eta' + g(t)\eta + h(t) = 0 \ph (\eta' + g_1(t) \eta + h_1(t) = 0)$ on $[t_1,+\infty)$ is also a solution of the inequality (2,3)  ((2.4)). Therefore if  $f(t) \ge 0 \phantom{a} (f_1(t) \ge 0), \phantom{a} t\ge t_0$, then the inequality  (2.3)   ((2.4)) has a solution on  $[t_0,+\infty)$,   satisfying any initial condition.

{\bf Theorem 2.1}.   {\it Let $y_1(t)$  be a solution of Eq. $(2.2)$  on  $[t_0,+\infty)$, and let  $\eta_0(t)$,  $\eta_1(t)$  be  solutions of Ineq. (2.3) and  Ineq. (2.4) respectively with  $\eta_k(t_1) \ge y_0(t_1), \phantom{a} k=0,1$. Let $f(t) \ge 0,$
$$
\lambda - y_1(t) + \int\limits_{t_1}^t \exp\biggl\{\int\limits_{t_1}^\tau[f(\xi)(\eta_0(\xi) + \eta_1(\xi)) + g(\xi)]\biggr\}\times \phantom{aaaaaaaaaaaaaaaaaaaaaaaaaaaaaa}
$$
$$
\phantom{aaaaaaaaa}\times [(f_1(t) - f(t)) y_1^2(t) + (g_1(t) - g(t)) y_1(t) + h_1(t) - h(t)]d\tau \ge 0, \phantom{aaa} t\ge t_0,
$$
for some $\lambda\in [y_1(t_1), \eta_0(t_1)]$. Then for each  $y_{(0)} \ge y_1(t_0)$
Eq. (2.1) has the solution  $y_0(t)$  on  $[t_0,+\infty)$,
satisfying the initial condition $y_0(t_1) = y_{(0)}$,   moreover   $y_0(t) \ge y_1(t), \phantom{a} t\ge t_0.$}

See the proof in  [13].

{\bf Definition 2.1.} {\it A solution of Eq. (2.1) is called $t_1$-regular if it exists on $[t_1,+\infty)$.}

{\bf Definition 2.2.} {\it A $t_1$-regular solution $y_0(t)$ of Eq. (2.1) is called $t_1$-normal if there exists a $\delta$-neighborhood $U_\delta(y_0(t_1)) \equiv (y_0(t_1) - \delta, y_0(y_1) + \delta) \ph (\delta > 0)$ of $y_0(t_1)$ such that every solution $y(t)$ of Eq. (2.1) with $y(t_1) \in U_\delta(y_0(t_1))$ is also $t_1$-regular, otherwise $y(t)$ is called $t_1$-extremal.}

Let $y(t)$ be a $t_1$-regular solution of Eq. (2.1). We can interpret $y(t)$ as a solution of the linear equation
$$
y' + G(t) y + h(t) = 0, \phh t\ge t_1,
$$
where $G(t) \equiv f(t) y(t) + g(t), \ph t\ge t_1$. Then by Cauchy formula we hve
$$
y(t) \equiv \exp\biggl\{-\il{i_1}{t}[f(\tau)y(\tau) + g(\tau)]d \tau\biggr\}\biggl[y(t_1) - \il{t_1}{t}\exp\biggl\{\il{t_1}{\tau}g(s)d s\biggr\} h(\tau)\phi_0(\tau)d \tau\biggr], \ph t \ge t_1,
$$
where $\phi_0(t) \equiv \exp\biggl\{\il{t_1}{t}f(\tau) y(\tau) d \tau\biggr\}, \ph t \ge t_1.$ From here it follows
$$
y(t)\phi_0(t) = y(t_1)\exp\biggl\{-\il{t_1}{t} g(\tau) d \tau\biggr\} - \il{t_1}{t}\exp\biggl\{-\il{\tau}{t}g(s) d s\biggr\} h(\tau)\phi_0(\tau) d \tau,  \phh t \ge t_1. \eqno (2.5)
$$

{\bf Lemma 2.1.} {\it Let $y(t)$ be a $t_1$-regular solution of Eq. (2.1) and let $f(t) \ge 0, \ph t \ge t_1.$ Then
$$
\il{t_1}{t} f(\tau) y(\tau) d \tau \le y(t_1)\il{t_1}{t} f(\tau) \exp\biggl\{-\il{t_1}{\tau}g(s) d s\biggr\}d \tau - \il{t_1}{t}f(\tau) d\tau\il{t_1}{\tau}\exp\biggl\{-\il{\xi}{\tau}g(s) d s\biggr\} h(\xi) d \xi,
$$
$t \ge t_1.$}

See the proof in [11].

{\bf Lemma 2.2.} {\it Let the following conditions be satisfied.
$$
\ilp{t_0} g(\tau) d \tau = +\infty, \ph \il{t_0}{t}\exp\biggl\{-\il{\tau}{t} g(s) d s\biggr\} |h(\tau)| d\tau \ph is \ph bounden \ph on \ph  [t_0,+\infty).
$$
Then for every continuous function $\phi(t)$, vanishing at $+\infty$, the relation
$$
\lim\limits_{t \to +\infty}\il{t_0}{t} \exp\biggl\{-\il{\tau}{t} g(s) d(s)\biggr\}|h(\tau)|\phi(\tau) d \tau = 0 \ph is \ph valid.
$$
}

See the proof in [11].

{\bf Lemma 2.3}. {\it Let $M=U diag\{m_1, ..., m_n\}U^*$, where $U$ is an unitary matrix, $U^*$ is the transpose conjugate of $U, \ph m_l\in\mathbb{H}, \ph l=\overline{1,n}.$ Then
$$
||M|| \le \exp\Bigl\{\max\{\re m_1, ..., \re m_n\}\Bigr\}.
$$
}
 Proof. Obviously $\exp\{M\} = U diag \{ e^{m_1}, ..., e^{m_n}\} U^*$ and $|e^{m_l}| = e^{\re m_l}, \ph l=\overline{1,n}.$ Moreover $||U|| = ||U^*|| = 1$. Hence, $||\exp\{M\}|| \le ||diag\{\{e^{m_1}, ..., e^{m_n}\}|| \le \max\limits_{l=\overline{1,n}}|e^{m_l}|= \linebreak  = \exp\Bigl\{\max\{\re m_1, ..., \re m_n\}\Bigr\}.$ The lemma is proved.

 Let $\mathfrak{a}_*(t)$ and  $\mathfrak{d}_*(t)$  real-valued continuous functions on $[t_0,+\infty)$ such that \linebreak $\re a_l(\tau;t) \le \il{t_0}{t}\mathfrak{a}_*(\tau) d \tau, \ph l = \overline{1,m} \ph \re  d_r(\tau;t) \le \il{t_0}{t}\mathfrak{d}_*(\tau)d\tau, \ph  r = \overline{1,n} \ph  t \ge \tau \ge t_0.$
Along with the system (1.1) consider the following scalar one
$$
\sist{\phi' = \mathfrak{a}_*(t) \phi + ||B(t)|| \psi,}{\phi' = ||C(t)|| \phi + \mathfrak{d}_*(t) \psi, \ph t\ge t_0.} \eqno (2.6)
$$

{\bf Lemma 2.4.} {\it Let the conditions $a)$ and $b)$ be satisfied. Then  the Lyapunov \linebreak (asymptotically) stability of the system (2.6) implies the Lyapunov (asymptotically) stability of the  system (1.1).}

Proof. Let $(\Phi(t), \Psi(t))$ be a solution of the system (1.1). We can interpret $\Phi(t)$ as a solution of the linear equation
$$
\Phi' = A(t) \Phi  + L(t), \phh t \ge t_0,
$$
where $L(t) \equiv B(t) \Psi(t), \phh t \ge t_0.$ Then by the condition $a)$ and the Cauchy formula we have
$$
\Phi(t) = \exp\biggl\{\il{t_0}{t}A(\tau) d \tau\biggr\}\biggl[\Phi(t_0) + \il{t_0}{t}\exp\biggl\{-\il{t_0}{\tau} A(s) d s\biggr\} L(\tau) d \tau\biggr], \phh t \ge t_0
$$
or
$$
\Phi(t) = \exp\biggl\{\il{t_0}{t}A(\tau) d \tau\biggr\}\biggl[\Phi(t_0) + \il{t_0}{t}\exp\biggl\{-\il{t_0}{\tau} A(s) d s\biggr\} B(\tau)\Psi(\tau) d \tau\biggr], \phh t \ge t_0. \eqno (2.7)
$$
By analogy for $\Psi(t)$ we can derive the equality
$$
\Psi(t) = \exp\biggl\{\il{t_0}{t}D(\tau) d \tau\biggr\}\biggl[\Psi(t_0) + \il{t_0}{t}\exp\biggl\{-\il{t_0}{\tau} D(s) d s\biggr\} C(\tau)\Phi(\tau) d \tau\biggr], \phh t \ge t_0. \eqno (2.8)
$$
Substitute in place of $\Psi(t)$ the right hand part of the last equality into (2.7). We  obtain
$$
\Phi(t) = F(t) + (K \Phi)(t), \phh t \ge t_0, \eqno (2.9)
$$
where
$$
F(t) \equiv \exp\biggl\{\il{t_0}{t} A(\tau) d \tau\biggr\}\Phi(t_0) + \phantom{aaaaaaaaaaaaaaaaaaaaaaaaaaaaaaaaaaaaaaaaaaaaaaaaaaaaaaaaaaa}
$$
$$
\phantom{aaaaaaaaaaa}+\expb{\il{t_0}{t}A(\tau) d\tau}\il{t_0}{t}\exp\biggl\{-\il{t_0}{\tau} A(s) ds\biggr\}B(\tau)\expb{\il{t_0}{\tau} D(s) d s}d\tau \Psi(t_0),
$$
$$
K(\Phi)(t) \equiv  \expb{\il{t_0}{t}A(\tau) d\tau}\il{t_0}{t}\exp\biggl\{-\il{t_0}{\tau} A(s) ds\biggr\}
B(\tau)  \expb{\il{t_0}{\tau}D(s)ds} d\tau \times \phantom{aaaaaaaaaaaaaaaaaaaaaaaaaaaaaaaa}
$$
$$
\phantom{aaaaaaaaaaaaaaaaaaaaaaaaaaaaaaaaaaaaaaa}\times \il{t_0}{\tau}\expb{-\il{t_0}{s}D(\xi) d \xi} C(s) \Phi(s) d s, \ph t\ge t_0.
$$
By analogy substituting in place of $\Phi(t)$  the right hand part of the equality (2.7) into (2.8) we arrive at the equality
$$
\Psi(t) = G(t) + (L\Psi)(t), \phh t \ge t_0, \eqno (2.10)
$$
where
$$
G(t) \equiv \exp\biggl\{\il{t_0}{t} D(\tau) d \tau\biggr\}\Psi(t_0) + \phantom{aaaaaaaaaaaaaaaaaaaaaaaaaaaaaaaaaaaaaaaaaaaaaaaaaaaaaaaaaaa}
$$
$$
\phantom{aaaaaaaaaaa}+\expb{\il{t_0}{t}D(\tau) d\tau}\il{t_0}{t}\exp\biggl\{-\il{t_0}{\tau} D(s) ds\biggr\}C(\tau)\expb{\il{t_0}{\tau} A(s) d s}d\tau \Phi(t_0),
$$
$$
L(\Psi)(t) \equiv  \expb{\il{t_0}{t}D(\tau) d\tau}\il{t_0}{t}\exp\biggl\{-\il{t_0}{\tau} D(s) ds\biggr\}
C(\tau)  \expb{\il{t_0}{\tau}A(s)ds} d\tau \times \phantom{aaaaaaaaaaaaaaaaaaaaaaaaaaaaaaaa}
$$
$$
\phantom{aaaaaaaaaaaaaaaaaaaaaaaaaaaaaaaaaaaaaaa}\times \il{t_0}{\tau}\expb{-\il{t_0}{s}A(\xi) d \xi} C(s) \Psi(s) d s, \ph t\ge t_0.
$$

Let $(\phi_0(t), \psi_0(t))$ be a solution of the system (2.6). By (2.9) and (2.10) we have respectively
$$
\phi_0(t) = F_0(t) + \il{t_0}{t} K_0(t,\xi) \phi_0(\xi) d \xi, \phh t \ge t_0, \eqno (2.11)
$$
$$
\psi_0(t) = G_0(t) + \il{t_0}{t} L_0(t,\xi) \psi_0(\xi) d \xi, \phh t \ge t_0, \eqno (2.12)
$$
where
$$
F_0(t) \equiv \phi_0(t_0)\exp\biggl\{\il{t_0}{t} \mathfrak{a}_*(\tau) d \tau\biggr\} + \psi_0(t_0)\il{t_0}{t}\exp\biggl\{\il{t_0}{\tau}\mathfrak{d}_*(s) ds + \il{\tau}{t}  \mathfrak{a}_*(s) d s\biggr\} ||B(\tau)|| d \tau,
$$
$$
t\ge t_0, \ph K_0(t,\xi) \equiv ||C(\xi)||\exp\biggl\{-\il{t_0}{\xi} \mathfrak{d}_*(s) d s\biggr\}\il{\xi}{t}\exp\biggl\{\il{\tau}{t} \mathfrak{a}_*(s) d s\biggr\} ||B(\tau)|| d \tau, \phh t \ge \xi \ge t_0,
$$
$$
G_0(t) \equiv \psi_0(t_0)\exp\biggl\{\il{t_0}{t} \mathfrak{d}_*(\tau) d \tau\biggr\} + \phi_0(t_0)\il{t_0}{t}\exp\biggl\{\il{t_0}{\tau} \mathfrak{a}_*(s) ds + \il{\tau}{t} \mathfrak{d}_*(s) d s\biggr\} ||C(\tau)|| d \tau,
$$
$$
t \ge t_0, \ph L_0(t,\xi) \equiv ||B(\xi)||\exp\biggl\{-\il{t_0}{\xi} \mathfrak{a}_*(s) d s\biggr\}\il{\xi}{t}\exp\biggl\{\il{\tau}{t} \mathfrak{d}_* d(s) d s\biggr\} ||C(\tau)|| d \tau, \phh t \ge \xi \ge t_0.
$$
By (2.9) and (2.10) we can represent $\phi(t)$ and $\psi(t)$ respectively via the following series expansion
$$
\Phi(t) = F(t) + (K F)(t)  +  (K^2 F)(t)  + ... , \phh t \ge t_0, \eqno (2.13)
$$
$$
\Psi(t) = G(t) + (L G)(t)  +  (L^2 G)(t)  + ... , \phh t \ge t_0, \eqno (2.14)
$$

Assume $\phi_0(t_0) = ||\Phi(t_0)||, \ph \psi_0(t_0) = ||\Psi(t_0)||$. Then by Lemma 2.3   and by (2.11) and (2.12) from  $b)$, (2.13) and (2.14) we obtain respectively:
$$
||\Phi(t)|| \le F_0(t) + \il{t_0}{t} K_0(t,\xi) F_0(\xi) d \xi + \il{t_0}{t} K_0(t,\xi)d\xi \il{t_0}{\xi} K_0(\xi,\zeta) F_0(\zeta) d \zeta + ... = \phi_0(t),  \eqno (2.15)
$$
$t \ge t_0$,
$$
||\Psi(t)|| \le  G_0(t) + \il{t_0}{t} L_0(t,\xi) G_0(\xi) d \xi + \il{t_0}{t} L_0(t,\xi)d\xi \il{t_0}{\xi} L_0(\xi,\zeta) G_0(\zeta) d \zeta + ... =\psi_0(t),  \eqno (2.16)
$$
$t \ge t_0.$ Assume the system (2.6) is Lyapunov (asymptotically) stable Then the estimates (2.15) and (2.16)  imply that the system (1.1) is also Lyapunov (asymptotically) stable. The lemma is proved.

{\bf 3. Main result.}  Set $E(t) \equiv \mathfrak{a}_*(t) - \mathfrak{d}_*(t), \ph t\ge t_0.$

{\bf Theorem 3.1.} {\it Let the conditions $a)$, $b)$ and the  conditions

\noindent
1) $\supl\il{t_0}{t}\expb{\il{\tau}{t} \mathfrak{d}_*(s) d s}||C(\tau)|| d \tau < +\infty$;

\noindent
2) $\supl\il{t_0}{t}\biggl[\mathfrak{a}_*(\tau) + ||B(\tau)||\il{t_0}{\tau}\expb{-\il{\xi}{\tau} E(s) ds} ||C(\xi)|| d \xi\biggr] d\tau < + \infty$;

\noindent
$\biggl( 2') \ilp{t_0} E(\tau) d\tau = +\infty, \hskip 3pt \liml \il{t_0}{t}\biggl[\mathfrak{a}_*(\tau)  +||B(\tau)||\il{t_0}{\tau}\expb{-\il{\xi}{\tau} E(s) ds} ||C(\xi)|| d \xi\biggr] d\tau = - \infty. \biggr)$

\noindent
be satisfied.

\noindent
Then the system (1.1) is Lyapunov (asymptotically) stable.}

Proof. In virtue of Lemma 2.4 to prove the theorem it is enough to show that the system (2.5) is Lyapunov (asymptotically) stable.
Consider the Riccati equations
$$
y' + ||B(t)|| y^2 +  E(t) y - ||C(t)|| = 0, \phh t \ge t_0, \eqno (3.1)
$$
$$
y' + ||B(t)|| y^2 + E(t) y  = 0, \phh t \ge t_0. \phh\phh
$$
Applying Theorem 2.1 to these equations we conclude that for every $t_1\ge t_0$ and $\gamma \ge 0$ Eq. (3.1) has a solution $y(t)$ on $[t_1,+\infty)$ with $y(t_1) = \gamma$ and
$$
y(t) \ge 0, \phh t \ge t_1. \eqno (3.2)
$$
All solutions $y(t)$ of Eq. (3.1), existing on $[t_0,+\infty)$, are connected with solutions $(\phi(t), \psi(t))$ of the system (2.6) by relations (see [11]):
$$
\phi(t) = \phi(t_0)\expb{\il{t_0}{t}\Bigl[||B(\tau)|| y(\tau) + \mathfrak{a}_*(\tau)\Bigr]d\tau}, \ph \phi(t_0)\ne 0, \ph \psi(t)= y(t)\phi(t),  \eqno (3.3)
$$
$t \ge t_0.$ Let $y_0(t)$ be the solution of Eq. (3.1) with $y_0(t_0) = 0$. By already proven $y_0(t)$ exists on $[t_0,+\infty)$ and is non negative. Set $\phi_0(t) \equiv \expb{\il{t_0}{t}||B(\tau)||y_0(\tau) d \tau}, \ph t \ge t_0$. By (2.5) we have
$$
y_0(t)\phi_0(t) = \il{t_0}{t}\expb{-\il{\tau}{t} E(s) d s} ||C(\tau)|| \phi_0(\tau) d \tau, \phh t\ge t_0. \eqno (3.4)
$$
Let $(\phi(t), \psi(t))$ be the solution of the system (2.6) with $\phi(t_0) =1, \ph \psi(t_0) = 0$. Then by (3.3) we have
$$
\phi(t) =\expb{\il{t_0}{t}\Bigl[||B(\tau)|| y_0(\tau) + \mathfrak{a}_*(\tau)\Bigr]d\tau}, \phh \psi(t) = y_0(t)\phi(t), \phh t\ge t_0. \eqno (3.5)
$$
This together with (3.4) implies
$$
\psi(t) = \il{t_0}{t}\expb{-\il{\tau}{t} \mathfrak{d}_*(s) d s}||C(\tau)|| \phi(\tau) d\tau, \phh t\ge t_0. \eqno (3.6)
$$
In virtue of Lemma 2.1 from the first equality of (3.5) it follows
$$
0 < \phi(t) \le \expb{\il{t_0}{t}\Bigl[\mathfrak{a}_*(\tau) +||B(\tau)||\il{t_0}{\tau}\expb{-\il{\xi}{\tau} E(s) d s}||C(\xi)||d \xi\Bigr]d\tau}, \phh t \ge t_0. \eqno (3.7)
$$
Show that $(\phi(t), \psi(t))$ is bounded (vanish at $+\infty$). The condition 2) (the condition2') together with (3.7) implies that $\phi(t)$ is bounded (vanish at $+\infty$). From here and from (3.6) (by Lemma 2.2 from here and from (3.6)) it follows that $\psi(t)$ is also bounded (vanish at $+\infty$). Therefore $(\phi(t), \psi(t))$ is bounded (vanish at $+\infty$). Since $c(t)\not\equiv 0$ there exists $t_1 \ge t_0$ such that $y(t) > 0, \ph t\ge t_1$. Hence by the second equality of (3.5)
$$
\psi(t) > 0, \phh t\ge t_1. \eqno (3.8)
$$
Let $(\phi_1(t), \psi_1(t))$ be a solution of the system (2.6) such that $\phi_1(t_1) > 0, \ph \psi_1(t_1) > 0$ and
$
det
\begin{pmatrix}
\phi(t_1) & \psi(t_1)\\
\phi_1(t_1) & \psi_1(t_1)
\end{pmatrix} \ne 0
$.
Then $(\phi(t), \psi(t))$ and $(\phi_1(t), \psi_1(t))$ are linearly independent. Taking into account (3.5) and (3.8) we have
$$
\frac{\psi(t_1)}{\phi(t_1)} > 0, \phh \frac{\psi_1(t_1)}{\phi_1(t_1)} > 0. \eqno (3.9)
$$
Let $y_1(t)$ be the solution of Eq. (3.1) with $y_1(t_1) = \frac{\psi_1(t_1)}{\phi_1(t_1)}$. Then by the second inequality of (3.9) and by the already proven $y_1(t)$ is $t_1$-normal. By the first inequality of (3.9) and by the already proven $y(t)$ is also $t_1$-normal. Therefore (see [14])
$$
M\equiv \Supl{t\ge t_1}\biggl|\il{t_1}{t} ||B(\tau)||(y(\tau) - y_1(\tau))d\tau\biggr| < +\infty. \eqno (3.10)
$$
By (3.3) we have
$$
\phi_1(t) = \phi_1(t_1) \expb{\il{t_1}{t}\Bigl[||B(\tau)|| y_1(\tau) + \mathfrak{a}_*(\tau)\Bigr] d\tau} = \frac{\phi_1(t_1)}{\phi(t_1)}  \expb{\il{t_1}{t}\Bigl[||B(\tau)|| y(\tau) + \mathfrak{a}_*(\tau)\Bigr] d\tau}\times
$$
$$
\phantom{aaaaaaaaaaaaaaaaaaaaaaaaaaaaaaaaaaaaaa}\times \expb{\il{t_1}{t}||B(\tau)||\Bigl[y_1(\tau) - y(\tau)\Bigr]d\tau}, \phh t\ge t_1.
$$
This together with (3.10) implies
$$
0 < \phi_1(t) < \frac{\phi_1(t_1)}{\phi(t_1)}\exp\{M\}\phi(t), \phh t \ge t_1.
$$
Therefore $\phi_1(t)$ is bounded (vanish at $+\infty$). Hence since $(\phi(t), \psi(t))$ and $(\phi_1(t), \psi_1(t))$ are linearly independent to complete the proof of the theorem  it is enough to show that $\psi_1(t)$ is bounded (vanish at $+\infty$). Let $z(t)$ and $z_1(t)$ be the solutions of the Riccati equation
$$
z' + ||C(t)|| z^2 - E(t) z - ||B(t)|| = 0, \phh t\ge t_1
$$
with $z(t_1) = \frac{\phi(t_1)}{\psi(t_1)} > 0, \ph z_1(t_1) = \frac{\phi_1(t_1)}{\psi_1(t_1)} > 0.$ Then by already proven $z(t)$ and $z_1(t)$ are $t_1$-normal, and therefore (see [14])
$$
M_1\equiv \Supl{t\ge t_1}\biggl|\il{t_1}{t} ||C(\tau)||(z(\tau) - z_1(\tau))d\tau\biggr| < +\infty. \eqno (3.11)
$$
By (3.3) we have
$$
\psi_1(t) = \psi_1(t_1)\expb{\il{t_1}{t}\Bigl[||C(\tau)||z_1(\tau) + \mathfrak{d}_*(\tau)\Bigr] d\tau}= \phantom{aaaaaaaaaaaaaaaaaaaaaaaaaaaaaaaaaaaaa}
$$
$$
=\frac{\psi_1(t_1)}{\psi(t_1)}\psi(t_1)\expb{\il{t_1}{t}\Bigl[||C(\tau)||z(\tau) + \mathfrak{d}_*(\tau)\Bigr] d\tau}\expb{\il{t_1}{t}||C(\tau)||(z_1(\tau) - z(\tau))d\tau}, \ph t\ge t_1.
$$
This together with (3.11) implies that
$$
0 < \psi_1(t) \le \frac{\psi_1(t_1)}{\psi(t_1)}\exp\{M_1\}\psi(t), \phh t \ge t_1.
$$
Hence $\psi_1(t)$ is bounded (vanish at $+\infty$). The theorem is proved.

{\bf Remark 3.1.} {\it The conditions $a)$ and $b)$ in Theorem 3.1 can be replaced by the following ones

\noindent
$a') \ph A(t)\il{t_0}{t}A(\tau) d \tau = \il{t_0}{t}A(\tau) d\tau A(t), \ph   \ph D(t)\il{t_0}{t}D(\tau) d \tau = \il{t_0}{t}D(\tau) d\tau D(t), \ph t \ge t_0;$

\noindent
the matrices $\il{t_0}{t}A(\tau) d \tau$ and $\il{t_0}{t}D(\tau) d \tau$ are unitary equivalent to some diagonal \linebreak matrices,~ i. e.

\noindent
$b') \ph U_A(t) \il{t_0}{t}A(\tau) d \tau U_A(t)^* = diag\{a_1(t), ..., a_m(t)\},$

$\phantom{aaaaaaaaaaaaaaaaaaaaaaaaaaaaaaaaa} U_D(t) \il{t_0}{t}A(\tau) d \tau U_D(t)^* = diag\{d_1(t), ..., d_m(t)\}$,

\noindent
where $U_A(t)$ and $U_D(t)$ are some unitary matrix functions, $U_A(t)^*$ and $U_D(t)^*$ are their transpose conjugate respectively, $a_l(t), \ph d_r(t) \in \mathbb{H},\ph l=\overline{1,m}, \ph r = \overline{1,n}, \ph t \ge t_0, \ph \re a_1(t) = ... = \re a_m(t), \ph \re d_1(t)= ... = \re d_n(t), \ph t\ge t_0$.}

Let $p(t), \ph q(t)$ and $r(t)$ be complex-valued continuous functions on $[t_0,+\infty)$ and let $p(t)\ne 0, \ph t\ge t_0$. Consider the second order linear ordinary differential equation
$$
(p(t)\phi')' + q(t)\phi' + r(t)\phi = 0, \phh t \ge t_0. \eqno (3.12)
$$
Introducing a new variable $\psi = p(t)\phi'$ we reduce this equation to the system
$$
\sist{\phi' = \phh \frac{1}{p(t)} \psi,}{\psi' = - r(t)\phi - \frac{q(t)}{p(t)} \psi, \ph t \ge t_0.} \eqno (3.13)
$$

{\bf Definition 3.1.} {\it Eq. (3.12) is called Lyapunov (asymptotically) stable if the system (3.13) is Lyapunov (asymptotically) stable.}

From Theorem 3.1 we immediately get

{\bf Corollary 3.1.} {\it Let the functions
$$
I_1(t)\equiv \il{t_0}{t}\expb{-\il{\tau}{t}\re\frac{q(s)}{p(s)}d s}|r(\tau)| d\tau, \phantom{aaaaaaaaaaaaaaaaaaaaaaaaaaaaaaaaaa}
$$
$$
\phantom{aaaaaaaaaaaaaaaaaaa}  I_2(t)\equiv \supl\il{t_0}{t}\frac{d\tau}{|p(\tau)|}\il{t_0}{\tau}\expb{-\il{\xi}{\tau}\re \frac{q(s)}{p(s)}}|r(\xi)|d\xi, \phh t\ge t_0
$$
 be bounded. Then Eq. (3.12) is Lyapunov stable.}

\phantom{aaaaaaaaaaaaaaaaaaaaaaaaaaaaaaaaaaaaaaaaaaaaaaaaaaaaaaaaaaaaaaaaaaaaaaa}$\blacksquare$

{\bf Remark 3.2.} {\it In the case $p(t) > 0, \ph r(t) \le 0$ and $q(t)$ is real-valued the condition of Corollary 3.1 is also necessary for Lyapunov stable of Eq. (3.12) (see [12]). In this sense the conditions 1) and 2) of Theorem 3.1 are sharp.}

{\bf Remark 3.3.} {\it It is not difficult to verify that in the case $p(t) > 0, \ph r(t) \le 0$ and $q(t)$ is real-valued the condition 2') of Theorem 3.1 for the system (3.13) is not satisfiable. On the other hand using Theorem 2.1 to the pair of equations
$$
y' + \frac{1}{p(t)} y^2 + \frac{q(t)}{p(t)} y + r(t) = 0, \phh t\ge t_0,
$$
$$
y' + \frac{1}{p(t)} y^2 + \frac{q(t)}{p(t)} y  = 0, \phh t\ge t_0 \phh \phh
$$
one can easily show that in this case Eq. (3.12) cannot be asymptotically stable (it has a positive and non decreasing solution). In This sence the condition 2') of Theorem 3.1 is sharp.}

{\bf Example 3.1.} {\it Consider the system
$$
\sist{\phi' = \nu(t)\phi + \frac{\mu(t)}{t \ln^2 t}\psi,}{\psi'= \mu(t)\phi + (\nu(t) - 1)\psi, \ph t\ge e,} \eqno (3.14)
$$
where $\nu(t)$ and $\mu(t)$ are some real-valued continuous functions on $[e,+\infty)$ and $\mu(t)$ is bounded on $[e,+\infty)$. Assume $\il{e}{t}\nu(\tau) d\tau$ is upper bounded on $[e,+\infty)$ ($\ilp{e}\nu(\tau) d \tau = -\infty$ and $\il{e}{t}(\varepsilon - \nu(\tau)) d\tau$ is upper bounded on $[e,+\infty)$ for some $\varepsilon \in (0,1)$). Then it is not difficult to verify that the conditions 1) and 2) (2')) of Theorem 3.1  for the system (3.14) are satisfied. Therefore under the indicated restrictions the system (3.14) is Lyapunov (asymptotically) stable. Since at least one of the integrals $\ilp{e}|\nu(\tau)|d\tau, \ph  \ilp{e}|\nu(\tau) -1|d\tau$ diverges to $+\infty$ the application of the estimates of Lyapunov and Yu. S. Bogdanov ([4], p. 133) to the system (3.14) gives no result. Let $A(t)$ be the matrix of the coefficients of the system (3.14). Then
$$
\gamma_\pm(t) \equiv \frac{2\nu(t) - 1 \pm \sqrt{1 + \frac{4\mu^2(t)}{t\ln^2 t}}}{2}, \phh t \ge e
$$
are its eigenvalues. Therefore if $\Supl{t \ge e} \nu(t) \ge0$, then the application of the freezing method ([4], p. 139, Theorem 4.6.4) to the system (3.14) gives no result. Let us now  discuss the applicability of estimates of solutions via logarithmic norms$\gamma_I(t), \ph \gamma_{II}(t)$ and $\gamma{III}(t)$ of S. M. Lozinski ([4], pp. 135, 136). From the Lozinski's theorem ([4], p. 137) it follows that if one of the integrals $\il{e}{t}\gamma_i(\tau) d\tau, \ph i= I, II, III$ is upper bounded then the corresponding linear system is Lyapunov stable. For the system (3.14) we have
$$
\gamma_I(t) \ge \nu(t) + |\mu(t)|, \phh t \ge e.
$$
Therefore if $\Supl{t\ge e}\il{e}{t}(\nu(\tau) + |\mu(\tau)|) d\tau = +\infty$, then the application of $\gamma_I(t)$ to the system (3.14) gives no result. If $|\mu(t)| \ge \frac{e}{e - 1}, \ph t \ge e$ then the logarithmic norm $\gamma_{II}(t)$ of the system (3.14) satisfies to the inequality
$$
\gamma_{II}(t) \ge 1 + \nu(t), \phh t \ge e.
$$
Hence if $\Supl{t \ge e}\il{e}{t}(1 + \nu(\tau)) d \tau = +\infty$, then the application of $\gamma_{II}(t)$ to the system  (3.14) gives no result. Finally the logarithmic norm $\gamma_{III}(t)$ for the system (3.14) is
$$
\gamma_{III}(t) = \frac{2\nu(t) - 1 + \sqrt{1 + 4\mu^2(t)(1 + \frac{1}{t\ln^2 t})}}{2}, \phh t \ge e.
$$
Therefore if $\ilp{e}\mu^2(\tau) d \tau = +\infty$ and $\il{e}{t}\nu(\tau) d \tau$ is bounded from below then \linebreak $\ilp{e}\gamma_{III}(\tau) d\tau =~ + \infty$ and, hence the application of $\gamma_{III}(t)$ to the system (3.14) gives also no result. Thus if $\il{e}{t}\nu(\tau) d\tau$ is bounded and  $|\mu(t)| \ge \frac{e}{e - 1}, \ph t \ge e$ then none of the logarithmic norms $\gamma_I(t), \ph \gamma_{II}(t)$ and $\gamma_{III}(t)$ is applicable to the system (3.14).}

{\bf Example 3.2.} {\it Consider the system
$$
\sist{\phi' = (\lambda - C\sin t) \phi + \mu_1\psi,}{\psi' = \phh \mu_2 \phi \phh +\lambda_2 \psi, \phh t \ge 0,} \eqno (3.15)
$$
where $\lambda_k, \ph \mu_k, \ph k=1,2, \ph C$ are some real constants, $\mu_k > 0, \ph k=1,2, \ph C > 0$. It is not difficult to verify that under the restrictions
$$
\lambda_k < 0, \ph k=1,2, \ph \lambda_1 - \lambda_2 > 0, \ph \lambda_1 + \frac{\mu_1\mu_2}{\lambda_1 - \lambda_2} \le 0 \ph (< 0) \eqno (3.16)
$$
the conditions 1) and 2) (2')) of Theorem 3.1 for the system (3.15) are satisfied. Therefore under these restrictions the system (3.15) is Lyapunov (asymptotically) stable. Since $\ilp{0}|\lambda_1 - C\sin t|d t = +\infty$ the application of estimates of Lyapunov and  Yu. S. Bogdanov gives no result. Let $A_1(t)$ be the matrix of coefficients of (3.15). Then
$$
\gamma(t) \equiv \frac{\lambda_1 - \lambda_2 - C\sin t + \sqrt{(\lambda_2 - \lambda_1 + C\sin t)^2 + 4\mu_1\mu_2}}{2}, \ph t \ge 0
$$
is its greatest eigenvalue. Hence, if $C\ge |\lambda_1 + \lambda_2|$ then $\Supl{t \ge 0} \gamma(t) \ge 0$, and in this case the freezing method is not applicable to the system (3.15). For (3.15) we have the following logarithmic norms:
$$
\gamma_I(t) = \max\{\lambda_1 + \mu_1 - C\sin t, \lambda_2 + \mu_2\};
$$
$$
\gamma_{II}(t) = \max\{\lambda_1 + \mu_2 - C\sin t, \lambda_2 + \mu_1\};
$$
$$
\gamma_{III}(t) = \frac{\lambda_1 + \lambda_2 - C\sin t + \sqrt{(\lambda_2 - \lambda_1 + C \sin t)^2 + (\mu_1 + \mu_2)^2}}{2}, \phh t \ge 0.
$$
The set of parameters $\mu_k, \ph \lambda_k$ for which at least one of norms   $\gamma_I(t), \ph \gamma_{II}(t)$
is applicable to (3.15) does not include the set, defined by (3.16).
For example for $\lambda_1 = -1, \ph \lambda_2 =~ -\frac{3}{2}, \linebreak \mu_1 = 2, \ph \mu_2 = \frac{1}{5}$ the application of  $\gamma_1(t)$ and $\gamma_{II}(t)$ to (3.15) gives no result, whereas for these values of $\lambda_k, \ph \mu_k, \ph k=1,2$ the conditions (3.16) are satisfied. It is not difficult to verify that for all enough large (with respect to $\lambda_k|, \ph \mu_k, \ph k=1,2$) $C$ the equality $\ilp{0}\gamma_{III}(t) d t=+\infty$ is satisfied. Therefore for all enough large $C$ the application of $\gamma_{III}(t)$ to the system (3.15) gives no result.}

\vskip 20pt

\centerline{\bf References}

\vskip 20pt

\noindent
1. L. S. Pontriagin, Obyknovennye differential'nye uravneniya (Ordinary differential  \linebreak \phantom{aa} equations) Moskaw, Nauka, 1974.

\noindent
2. L. Cezary, Asymptotic Behavior and Stability Problems in Ordinary Differential \linebreak \phantom{aa} Equations, Berlin, 1959.

\noindent
3. N. W. Mac Lachlan, Theory and application of Mathieu Functions, Oxford, Clarendon \linebreak \phantom{aa} Press, 1947.

\noindent
4. L. Y. Adrianoba,  Introduction to the theory of linear systems of differential equations. \linebreak \phantom{aa}
S. Peterburg, Publishers of St. Petersburg University, 1992.

\noindent
5. V. A. Yakubovich, V. M. Starzhinsky,  Linear differential equations with periodic \linebreak \phantom{aa}  coefficients and their applications. Moscow, ''Nauka'', 1972.

\noindent
6. R. Bellman,  Stability theory of differential equations.
 Moscow, Foreign Literature \linebreak \phantom{aa} Publishers, 1954.

\noindent
7. V. I. Burdina,  On boundedness of solutions of systems of differential equations, Dokl.\linebreak \phantom{aa}
AN SSSR, 93:4 (1953), 603–606.

\noindent
8. I. M. Sobol. Study of the asymptotic behaviour of the solutions of the linear \linebreak \phantom{aa} second order differential equations wit the aid of polar coordinates. "Matematicheskij \linebreak \phantom{aa} sbornik", vol. 28 (70), N$^\circ$ 3, 1951, pp. 707 - 714.

\noindent
9. M. V. Fedoriuk. Asymptotic methods for linear ordinary differential equations. \linebreak \phantom{aa} Moskow, ''Nauka'', 1983.

\noindent
10.  Ph. Hartman, Ordinary differential Equations.
Second Edition, SIAM, 2002.

\noindent
11. G. A. Grigorian, On the Stability of Systems of Two First - Order Linear Ordinary \linebreak \phantom{aa} Differential Equations, Differ. Uravn., 2015, vol. 51, no. 3, pp. 283 - 292.

\noindent
12. G. A. Grigoryan, Stability Criterion for Systems of Two First-Order Linear Ordinary \linebreak \phantom{aa}  Differential Equations.  Math. Notes, 103:6 (2018), 892–900.

\noindent
13. G. A. Grigorian.  On Two Comparison Tests for Second-Order Linear  Ordinary \linebreak \phantom{aa}   Differential Equations (Russian) Differ. Uravn. 47 (2011), no. 9, 1225 - 1240; translation \linebreak \phantom{aa}in Differ.
 Equ. 47 (2011), no. 9 1237 - 1252, 34C10.

\noindent
14. G. A. Grigorian, On some properties of solutions of the Riccati equation.
Izvestiya  NAS \linebreak\phantom{aa} of Armenia, vol. 42, $N^\circ$ 4, 2007, pp. 11 - 26.

\end{document}